\documentclass[10pt]{amsart}
\usepackage{latexsym, amsmath,amssymb}
\usepackage{amsmath}
\usepackage{amsfonts}
\usepackage{amssymb}
\usepackage{amscd}
\usepackage{amsthm}
\usepackage{amsbsy}
\usepackage{graphicx}
\usepackage{bm}

\renewcommand{\epsilon}{\varepsilon}
\newcommand{\newsection}[1]
{\subsection{#1}\setcounter{theorem}{0} \setcounter{equation}{0}
\par\noindent}

\newtheorem{theorem}{Theorem}

\newtheorem{lemma}[theorem]{Lemma}
\newtheorem{corr}[theorem]{Corollary}

\newtheorem{proposition}[theorem]{Proposition}
\newtheorem{deff}[theorem]{Definition}

\newcommand{\bth}{\begin{theorem}}
\newcommand{\ble}{\begin{lemma}}
\newcommand{\bcor}{\begin{corr}}

\newcommand{\bdeff}{\begin{deff}}

\newcommand{\bprop}{\begin{proposition}}
\newcommand{\ele}{\end{lemma}}
\newcommand{\ecor}{\end{corr}}
\newcommand{\edeff}{\end{deff}}

\newcommand{\eprop}{\end{proposition}}

\newcommand{\e}{\varepsilon}
\renewcommand{\l}{\lambda}

\renewcommand{\Pi}{\varPi}

\renewcommand{\epsilon}{\varepsilon}

\newcommand{\1}{{\rm 1\hspace*{-0.4ex}%
\rule{0.1ex}{1.52ex}\hspace*{0.2ex}}}

\newcommand{\E}{{\mathbf E}}

\begin{document}


\title[$L^4$ norms of typical eigenfunctions]
{Concerning the $L^4$ norms of typical eigenfunctions on compact surfaces}
\thanks{The authors were supported in part by the NSF}

\author{Christopher D. Sogge}
\author{Steve Zelditch}
\address{Johns Hopkins University, Baltimore, MD}
\address{Northwestern University, Evanston, IL}

\maketitle

\begin{abstract}
Let $(M,g)$ be a two-dimensional compact
boundaryless Riemannian manifold with Laplacian, $\Delta_g$.  If
$e_\lambda$ are the associated eigenfunctions of
$\sqrt{-\Delta_g}$ so that $-\Delta_g e_\lambda = \lambda^2
e_\lambda$, then it has been known for some time \cite{soggeest}
that $\|e_\lambda\|_{L^4(M)}\lesssim \lambda^{1/8}$, assuming that
$e_\lambda$ is normalized to have $L^2$-norm one.  This result is
sharp in the sense that it cannot be improved on the standard
sphere because of highest weight spherical harmonics of degree
$k$. On the other hand, we shall show that the average $L^4$ norm
of the standard basis for the space ${\mathcal H}_k$ of  spherical
harmonics of degree $k$ on $S^2$ merely grows like $(\log
k)^{1/4}$.  We also sketch a proof that the average of  $\sum_{j =
1}^{2k + 1} \|e_\lambda\|_{L^4}^4$ for a random orthonormal basis
of ${\mathcal H}_k$ is $O(1)$.
 We
are not able to determine the maximum of this quantity over all
orthonormal bases of ${\mathcal H}_k$  or for orthonormal bases of
eigenfunctions on  other Riemannian  manifolds. However, under the
assumption that the periodic geodesics in $(M,g)$ are of measure
zero, we are able to show that for {\it any} orthonormal basis of
eigenfunctions we have that
$\|e_{\lambda_{j_k}}\|_{L^4(M)}=o(\lambda_{j_k}^{1/8})$ for a
density one subsequence of eigenvalues $\lambda_{j_k}$.  This
assumption is generic and it is the one in the
Duistermaat-Gullemin theorem \cite{dg} which gave related
improvements for the error term in the sharp Weyl theorem. The
proof of our result uses a recent estimate of the first author
\cite{Sokakeya}
 that
gives a necessary and sufficient condition that
$\|e_\lambda\|_{L^4(M)}=o(\lambda^{1/8})$.
\end{abstract}


\newsection{Introduction}

The purpose of this note is to introduce a new problem on $L^p$
norms of eigenfunctions on compact Riemannian manifolds $(M,g)$.
We   prove some initial
  results on the problem,  and also include
some conjectures and heuristic remarks.

The problem, roughly speaking, is to determine the asymptotic
average of the  $L^4$ norms $||e_{\lambda}||_4$ of the elements of
an orthonormal basis of eigenfunctions
$$-\Delta_g e_\lambda = \lambda^2 e_\lambda$$
of the associated Laplace-Beltrami operator. In practice it is
simpler to consider the fourth power Weyl sums,
\begin{equation} \label{WEYL} \frac{1}{N(\lambda)} \sum_{j: \lambda_j \leq \lambda}
||e_{\lambda_j}||_4^4 \end{equation} where
$$N(\lambda)=\#\{\lambda_j\le \lambda\},$$
is the Weyl counting function. The asymptotics of (\ref{WEYL})
depend on the entire orthonormal basis and, as will be seen below,
can behave quite differently from the behavior of individual
eigenfunctions in the basis.

Before stating our results, let us recall the results on $L^p$
norms of individual eigenfunctions. In 1988, one of us showed in
\cite{soggeest} that for $2<q\le \infty$ and
\begin{multline}\label{1.1}\sigma(q)=\max \bigl(\, 2(1/2-1/q)-1/2, \tfrac12(1/2-1/q)\, \bigr)
\\
=
\begin{cases}
2(1/2-1/q)-1/2, \quad q\ge 6
\\
\tfrac12(1/2-1/q), \quad 2< q\le 6,
\end{cases}
\end{multline}
we have
\begin{equation}\label{1.2}
\|e_\lambda\|_{L^q(M)}\lesssim \lambda^{\sigma(q)},
\end{equation}
assuming as we shall do throughout that the eigenfunctions
are $L^2$-normalized so that
$$\|e_\lambda\|_{L^2(M)}=1,$$
where the norms are taken with respect to the volume element, $dV$.
This result is sharp since certain spherical harmonics on the
sphere, $S^2$, with the round metric saturate the estimate \eqref{1.2}.
Specifically, when $q\ge6$, $L^2$-normalized zonal functions,
$Z_k$ satisfy
$$\|Z_k\|_{L^q(S^2)}\approx k^{2(1/2-1/q)-1/2}, \quad q\ge6,$$
while the $L^2$-normalized highest weight spherical harmonics,
$Q_k=c_k(x_1+ix_2)^k$ satisfy
$$\|Q_k\|_{L^q(S^2)}\approx k^{\tfrac12(1/2-1/q)}, \quad q\ge 2.$$
Both are eigenfunctions of the standard Laplacian on $S^2$ with
eigenvalue $\lambda^2=k(k+1)$ in the above notation. Also, we are
taking $S^2$ to be $\{(x_1,x_2,x_3): \, x_1^2+x_2^2+x_3^2=1\}$, so
that, as $k\to \infty$, the $Q_k$ become highly concentrated on
the equator where $x_3=0$. The orthonormal basis of joint
eigenfunctions of $\Delta_g$ and of $x_3$-axis rotations are
generally denoted by $Y^k_m, m = - k, \dots, k$; in particular,
$Z_k  = Y^k_0$ and $Q_k = Y_k^k$.

Even though \eqref{1.1} cannot be improved on the sphere, it is
thought that for generic manifolds one has at least
\begin{equation}\label{1.3}
\|e_\lambda\|_{L^q(M)}=o(\lambda^{\sigma(q)}),
\quad \text{as } \, \, \lambda\to \infty,
\end{equation}
for a given $q>2$. In \cite{soggezelditch} we showed that for
generic $(M,g)$ this is true for $q>6$ (and also corresponding
results for higher dimensions).  This just followed from showing
that under a certain generic condition on $(M,g)$ one can improve
the $L^\infty$ estimate in \eqref{1.2} to be $\|e_\lambda\|_\infty
=o(\lambda^{1/2})$, which implies \eqref{1.3} for all $q>6$ by
interpolating with \eqref{1.2} for $q=6$.  The results in
\cite{soggezelditch} were recently improved in \cite{stz}. The key
point was to show that the bound $\|e_\lambda\|_\infty =
O(\lambda^{1/2})$ can only be obtained on $(M, g)$ possessing a
``peak  point" or ``pole" $z_0$ with the property that a positive
measure of directions in $S^*_{z_0} M$ exponentiate to geodesic
loops which return to $z_0$ at some time. This behavior  occurs at
poles of a surface of revolution, since all meridians are closed
geodesics through the pole,  and in particular explains why the
sup norm bounds are attained by zonal functions$Z_k$ on the round
sphere (see \S 3 below).

Even though there are satisfactory results concerning
\eqref{1.3} for relatively large exponents $q>6$, much
less is known for relatively small exponents $2<q<6$.
In this case, it is thought that the enemy for \eqref{1.3}
is maximal concentration along periodic geodesics,
as occurs for the highest weight spherical harmonics.
Using the formula for the $Q_k$ one checks that
they have $L^2$-mass bounded below on shrinking
$k^{-1/2}$
neighborhoods of the equator $\gamma=
\{(x_1,x_2,x_3)\in S^2: \, x_3=0\}$.  In \cite{Sokakeya},
(following an earlier result in \cite{bourgainef}), the first
author proved that for $2<q<6$, \eqref{1.3} is valid if and
only if this type of concentration does not occur.  Specifically,
a necessary and sufficient condition for \eqref{1.3} for this
range of exponents is that
\begin{equation}\label{4}
\sup_{\gamma\in \Pi} \int_{\text{dist}_g(\gamma, y)\le \lambda^{-1/2}}
|e_\lambda(y)|^2 \, dV = o(1),
\end{equation}
where $\Pi$ is the space of all unit-length geodesics in
$M$, and $\text{dist}_g(\, \cdot, \, \cdot\, )$ is the geodesic
distance associated to the metric $g$.

The goal of this paper is to show that even though on some
manifolds there are eigenfunctions
$e_\lambda$ having $L^4$-norms of maximal size $\approx \lambda^{1/8}$
as $\lambda\to \infty$, they are very sparse.  Our first result
of this type says that given any orthonormal basis $\{e_{\lambda_j}\}$
of eigenfunctions with eigenvalues $\lambda_1\le \lambda_2\le\dots$
on a two-dimensional compact Riemannian manifold $(M,g)$ with a zero
measure of periodic geodesics,
one can find a density one subsequence of eigenvalues, $\{\lambda_{j_k}\}$,
for which
\begin{equation}\label{1.5}
\|e_{\lambda_{j_k}}\|_{L^4(M)}=o(\lambda_{j_k}^{1/8}).
\end{equation}
By interpolation with the $L^6$-estimate  in \eqref{1.2}
and the trivial $L^2$-estimate, this implies
that we also have $\|e_{\lambda_{j_k}}\|_{L^q}=o(\lambda_{j_k}^{\sigma(q)})$ for
every $2<q<6$.  Presently, we do not how to prove the corresponding
results for $q\ge6$, or how to obtain any results like this for higher
dimensions $n\ge3$.  The condition that $(M,g)$ have a zero set of
periodic geodesics is generic and it is the assumption in the
 Duistermatt-Guillemin theorem \cite{dg}, which involved a similar $o$-improvement of the error term in the Weyl formula.

The assumptions that the periodic orbits are of measure zero of
course is not valid for the sphere. Nonetheless,  we can prove a
much stronger result for the standard basis $\{Y^k_m\}$ on $S^2$,
even though, as we pointed out before, this eigen-basis has
functions saturating \eqref{1.2} for each $2<q\le \infty$.

To be more specific, we recall that the Laplace-Beltrami operator on
$S^2\subset {\mathbb R}^3$ with the standard round metric
has eigenvalues $\lambda^2=k(k+1)$ repeating with multiplicity
$2k+1$, meaning that the corresponding eigenspace ${\mathcal H}_k$
of spherical harmonics of degree $k$ has this dimension.  If we use
longitudinal coordinates $\phi \in [0,\pi]$ and latitudinal ones
$\theta\in [0,2\pi]$ so that $S^2\ni x=(\sin \phi \cos\theta, \sin\phi \sin\theta, \cos\phi)$,
then in these coordinates the standard basis for ${\mathcal H}_k$ has
elements
\begin{equation}\label{1.6}
Y^k_m(\phi,\theta)=c_{m,k}P^m_k(\cos\phi)e^{im\theta}, \quad -k\le
m\le k,
\end{equation}
where $P^m_k$ are Legendre functions and $c_{m,k}$ are
$L^2$-normalizing constants.  When $m=0$, $Y^k_0$ is the zonal
function $Z_k$, and when $m=\pm k$ it is a highest weight
spherical harmonic of degree $k$.  For this basis, we shall show
that the average $L^4$-norm is of size $\approx (\log k)^{1/4}$,
as $k\to \infty$,  i.e.,
\begin{equation}\label{1.7}
\frac1{2k+1}\sum_{m=-k}^k\int_{S^2}|Y^k_m|^4\, dV \approx \log k,
\quad k\ge2,
\end{equation}
which of course is much stronger than \eqref{1.5} since it shows
that there must be a density one sequence of eigenfunctions among
this basis with $L^4$-norms growing logarithmically with respect
to the eigenvalues.  It seems somewhat paradoxical at first  that
\eqref{1.7} is valid for the standard basis on the sphere, while
the same basis  is the worst case for \eqref{1.2}. But this holds
because  the left side of \eqref{1.7} is a functional of  an
orthonormal basis rather than of individual eigenfunctions, and
most elements $Y^k_m$ have relatively small $L^4$ norms. It is
doubtful that the $\{Y^k_m\}$ maximize this functional among
orthonormal bases of spherical harmonics. In Section \ref{ONB} we
explain this further.

These observations raise the following \medskip

\noindent{\bf Problem}  Let $\dim M = 2$. For which $(M, g)$ (if
any) does there exist an orthonormal basis of eigenfunctions for
which there exists a positive density subsequence $e_{\l_{j_k}}$ so
that $||e_{\l_{j_k}}||_{L^4} = \Omega (\lambda_{j_k}^{1/8}).$? Or
is a result like \eqref{1.5} is  valid on any compact surface?
\medskip

We prove \eqref{1.7} by obtaining pointwise bounds for the $\ell^4(m)$ norms
of the basis elements of ${\mathcal H}_k$.  Specifically, we shall prove sharp
estimates for
$$\|Y^k_m(x)\|_{\ell^4(m)}=\left(\, \sum_{m=-k}^k|Y^k_m(x)|^4\, \right)^{1/4}.$$
By the inclusion $\ell^4\subset \ell^2$, this quantity is bounded by the
corresponding $\ell^2(m)$-norm.  On the sphere, the $\ell^2(m)$ norm is independent
of $x$, and, in fact,
\begin{equation}\label{8}
\|Y^k_m(x)\|_{\ell^2(m)}=\sqrt{(2k+1)/4\pi}.
\end{equation}
The $\ell^4(m)$ norm is of this order of of magnitude for points $x$ of distance
$O(1/k)$ from the poles where $\phi=0$ or $\pi$, but in order to obtain
\eqref{1.7} much better estimates are needed.  We shall obtain such an improvement,
which turns out to be sharp, by using \eqref{1.6} and well known asymptotics
for the kernel of the projection onto the spherical harmonics of degree
$k$, ${\mathcal H}_k$.  Thus, we are very much using here special properties
of $S^2$.   Our results can be thought of as a natural analog for $S^2$
of Zygmund's \cite{zygmund} theorem for the two-torus ${\mathbb T}^2$,
which says that the eigenfunctions of its Laplace-Beltrami operator have uniformly
bounded $L^4$-norms.  As we pointed out before, this is far from true on
$S^2$, but in an averaged sense it is almost true since the average $L^4$-norms
just grow like powers of logs of the eigenvalues.

For general Riemannian manifolds of dimension $n$, the local Weyl formula says
that if $N$ is large enough and fixed then
\begin{equation}\label{1.9}
\left(\, \sum_{|\lambda_j-\lambda|\le N}|e_{\lambda_j}(x)|^2\, \right)^{1/2}
\approx \lambda^{(n-1)/2}.
\end{equation}
It would be interesting to see to what extent there is an
improvement in the general case when one replaces this
$\ell^2$-norm by $\ell^q$ norms with $q>2$ and to what extent
results of this type perhaps depend on properties of the geodesic
flow  starting at $x$.  In a future work, we intend to carry out
the analysis for round spheres of dimension $n\ge3$ and certain
surfaces of revolution.  Understanding the case of general
manifolds and to what extent these results might depend on $x$
seems difficult at present.  On the other hand, by using estimates
like \eqref{1.2}, one can see that for most points $x\in M$, once
can improve on the trivial consequence of \eqref{1.9} that
$$\left(\, \sum_{|\lambda_j-\lambda|\le N}|e_{\lambda_j}(x)|^q\, \right)^{1/q}\lesssim
\lambda^{(n-1)/2}.$$
For instance, if $q=4$ and $n=2$, then using \eqref{1.2} and Tchebyschev's inequality
one sees that if $C<\infty$ is fixed then
$$\left|\, \left\{x\in M: \, \bigl(\, \sum_{|\lambda_j-\lambda|\le N}|e_{\lambda_j}(x)|^4\,
\bigr)^{1/4}\ge C\lambda^{1/2}\right\}\, \right|=O(\lambda^{-1/2}),$$
which, not surprisingly, is exactly the size of the sets on which the highest weight spherical
harmonics are concentrated.

This paper is organized as follows.  In the next section we shall present the
proof of \eqref{1.5}.  Then we shall turn our attention to the
sphere $S^2$ and prove the much stronger bounds \eqref{1.7} for $S^2$.

\newsection{$L^4$ norms of generic eigenfunctions}


In this section we shall establish \eqref{1.5}.  Specifically, we shall prove
the following

\begin{theorem}\label{theorem1.1}
Let $(M,g)$ be a two-dimensional compact Riemannian
manifold.  If $\Phi_t: S^*M\to S^*M$ is geodesic flow on
the cosphere bundle, assume that the set
\begin{equation}\label{2.0}
{\mathcal P}=
\{(x,\xi)\in S^*M: \Phi_t(x,\xi)
=(x,\xi), \, \, \text{some } t>0\}
\end{equation} has measure zero in $S^*M$ with respect to the volume
element.
Then if $e_{\lambda_j}$ is an orthonormal basis of eigenfunctions,
$-\Delta e_{\lambda_j}=\lambda_j^2$, with $\lambda_1\le \lambda_2\le\dots$
there is a subsequence of eigenvalues $\lambda_{j_k}$ satisfying
\begin{equation}\label{2.1}
\lim_{\lambda\to \infty}\frac{\#\{\, \lambda_{j_k}\le \lambda\}}{N(\lambda)}=1,
\end{equation}
so that
\begin{equation}\label{2.2}
\|e_{\lambda_{j_k}}\|_{L^4(M)}=o(\lambda^{1/8}_{j_k}).
\end{equation}

\end{theorem}

To prove this we shall use an estimate from \cite{Sokakeya} and arguments
from \cite{cdv} and \cite{zelditch}.  The estimate from \cite{Sokakeya} says
given $(M,g)$ as above
there is a uniform constant $C$ so
that if $-\Delta e_\l = \l^2 e_\l$ and $N=1,2,3,\dots$ then
\begin{multline*}
\int_M |e_\l (x)|^4 \, dV
\le CN^{-1/2}\l^{1/2}\|e_\l\|^4_{L^2(M)}
\\
+CN\l^{1/2}\|e_\l\|_{L^2(M)}\left[ \, \sup_{\gamma\in \Pi}
\int_{{\mathcal T}_{\l^{-1/2}}(\gamma)}|e_\l(x)|^2\, dV\, \right]
+C\|e_\l\|_{L^2(M)}^4.
\end{multline*}
Here $dV=dV_g$ is the volume element,  $\Pi$ is the space
of all unit-length geodesics, and
$${\mathcal T}_\varepsilon(\gamma)=\{y\in M: \, \text{dist}_g(y, \gamma)\le
\e \},$$
denotes an $\e$-tube about $\gamma$.  By optimizing the choice of $N$,
we see that the preceding inequality implies that
\begin{equation}\label{2.3}
\|e_\l\|_{L^4(M)}\le C\l^{1/8}\|e_\l\|_{L^2(M)}^{5/6}\, \sup_{\gamma\in\Pi}
\|e_\l\|^{1/6}_{L^2({\mathcal T}_{\l^{-1/2}})} + C\|e_\l\|_{L^2(M)}.
\end{equation}

In addition to this we require the following result which is a simple consequence
of the local Weyl law (see \cite{HoIV}).

\begin{lemma}\label{lemma1.1}
Let $(M,g)$ be a compact Riemannian manifold and let $A\in \Psi_{cl}^0(M)$ be
a classical pseudo-differential operator on $M$ of order zero.  Then
if $A_0$ is the principal symbol of $A$,
\begin{equation}\label{2.4}
\sum_{\lambda_j\le \lambda}\int_M |Ae_{\lambda_j}(x)|^2 dV
=(2\pi)^{-n}\lambda^n \int_{T^*B}|A_0(x,\xi)|dV d\xi
+O(\lambda^{n-1}).
\end{equation}
\end{lemma}

Here, $T^*B\subset T^*M$ is the ball bundle,  $\{(x,\xi): \sum_{jk}g^{jk}(x)\xi_j\xi_k\le 1\}$,
where $g^{jk}$ is the cometric, i.e., $(g^{jk}(x))^{-1}=(g_{jk}(x))$.  Note that \eqref{2.4} with
$A$ being the identity operator is the sharp Weyl formula (\cite{Av}, \cite{Le}, \cite{Ho1}), and the proof of the more general
case just follows from a straightforward modifications of that of this special case.

As a first step in the proof of Theorem~\ref{theorem1.1}, let us use some ideas from the proof of
the Duistermaat-Guillemin theorem \cite{dg} (see also \cite{ivrii}).  Given $(x,\xi)\in S^*M$
we define $L(x,\xi)$ for $(x,\xi)\in {\mathcal P}$ to be the minimal $t>0$ so that
$\Phi_t(x,\xi)=(x,\xi)$ and we define $L(x,\xi)$ to be $+\infty$ if $(x,\xi)\notin {\mathcal P}$,
where ${\mathcal P}$ is as in \eqref{2.0}.  Then $L(x,\xi)$ is clearly a lower semicontinuous
function on $S^*M$.    As a result, since we are assuming that ${\mathcal P}$ has measure
zero, it follows that for a given $T>0$
$${\mathcal P}_T=\{(x,\xi)\in S^*M: \, L(x,\xi)\le T\}$$
is a closed subset of $S^*M$ which is of measure zero since ${\mathcal P}_T\subset
{\mathcal P}$.  Therefore, given $\varepsilon>0$, we can find a
pseudodifferential operator $b\in \Psi^0_{cl}(M)$ whose
principal symbol satisfies $0\le b_0(x,\xi)\le1$, $b_0(x,\xi)=1$ for $(x,\xi)\in
{\mathcal N}({\mathcal P}_T),$ where ${\mathcal N}({\mathcal P}_T)$ is a neighborhood
of $\Pi$ in $S^*M$ and
$$\int_{B^*M}|b_0(x,\xi)|^2 dVd\xi <\varepsilon/3.$$
By Lemma~\ref{lemma1.1}, we conclude from this that
\begin{equation}\label{2.5}
\frac1{N(\lambda)}\sum_{\lambda_j\le \lambda}\int_M |be_{\lambda_j}|^2 \,
dV<\varepsilon/3+O_b(\lambda^{-1}),
\end{equation}
since we are assuming that $|T^*B|=1$.

If we let $B=Id-b\in \Psi^0_{cl}(M)$ then we claim that there is a uniform constant
$C$, which is independent of $\varepsilon$ and $T$ above so that
\begin{equation}\label{2.6}
\sup_{\gamma\in \Pi}\int_{{\mathcal T}_{\lambda_j^{-1/2}(\gamma)}} |Be_{\lambda_j}|^2
\, dV\le C/T+C'_{B,T}\lambda^{-1/2}.
\end{equation}
If $T$ is chosen  large enough so that $C/T<\varepsilon/3$, the preceding inequalities
imply that there is an $\lambda_0=\lambda_0(\varepsilon)$ so that
\begin{equation}\label{2.7}
\frac1{N(\lambda)}\sum_{\lambda_j\le \lambda}\sup_{\gamma\in \Pi}
\int_{{\mathcal T}_{\lambda_j^{-1/2}(\gamma)}}|e_{\lambda_j}(x)|^2 \, dV
<\varepsilon, \quad \text{if } \lambda>\lambda_0.
\end{equation}
As we shall see, this and \eqref{2.3} immediately yield Theorem~\ref{theorem1.1}.

Our main estimate \eqref{2.6} would follow from showing that there is a constant
$C$ as above so that
\begin{equation}\label{2.8}
\sup_{\gamma\in \Pi}\int_\gamma |Be_\lambda|^2\, ds\le CT^{-1}\lambda^{1/2}+C'_{T,B}.
\end{equation}
Here $ds$ is the geodesic arclength measure on $\gamma$.
Estimate \eqref{2.8} yields \eqref{2.9} due to the simple fact
that if $f\ge0$ then for any $\gamma_0\in \Pi$
$$\int_{{\mathcal T}_{\lambda^{-1/2}}(\gamma_0)} f\, dV
\le C\lambda^{-1/2}\sup_{\gamma\in\Pi}\int_\gamma f\, ds$$
for a uniform constant $C$ since ${\mathcal T}_{\lambda^{-1/2}}(\gamma_0)$
is a tube of width $\lambda^{-1/2}$ about $\gamma_0$.
\footnote{Note that, in ${\mathbb R^2}$,
the integral of $f\ge0$ over an $1\times \lambda^{-1/2}$ rectangle is dominated
by $\lambda^{-1/2}$ times the supremum of integrals over the line
segments in the rectangle that are parallel to the center segment, and
a similar argument works for the above tubes if one uses Fermi normal
coordinates about a geodesic which intersects $\gamma_0$ orthogonally.}
 Let
$\Pi_{cl}$ be the set of unit geodesics that are part of a
periodic geodesic. In \cite{Sokakeya} one of us showed that
$$\int_{\gamma}|e_\lambda|^2\, ds =o(\lambda^{1/2}), \quad
\text{if }\gamma\in \Pi\backslash \Pi_{cl},$$
which was an $o$-improvement of the restriction bounds in \cite{burq}.  The proof
of \eqref{2.8} is an adaptation of the one used to establish this result.

To prove \eqref{2.8}, let us fix a real-valued even function $\chi\in {\mathcal S}({\mathbb R})$
with $\chi(0)=1$ and $\Hat \chi(t)=0$, $|t|>1/4$, where $\Hat \chi$ denotes the
Fourier transform of $\chi$.  We then have that
$$\chi(T(P-\lambda))e_\lambda = e_\lambda$$
if $P=\sqrt{-\Delta_g}$.  Therefore, in order to prove \eqref{2.8}, it suffices to show that
\begin{equation}\label{2.9}
\int_\gamma |B\chi(T(P-\lambda))f|^2\, ds
\le CT^{-1}\lambda^{1/2}\|f\|^2_{L^2(M)}+C'_{T,B}\|f\|_{L^2(M)}^2,
\quad \gamma\in \Pi,
\end{equation}
where $C$ (but not $C'_{T,B}$) is a uniform constant independent
of $T$ and $B$.  We shall assume in what follows that $T$ is fixed but large,
in particular $T>10$.

Note that $\Pi$ is compact.  Therefore, in order to prove \eqref{2.9}, it suffices
to show that given $\gamma_0\in \Pi$ there is a neighborhood ${\mathcal N}(\gamma_0)$
 of $\gamma_0$ in $\Pi$ on which the analog of \eqref{2.9} holds with
constants independent of $\gamma\in {\mathcal N}(\gamma_0)$.  Different arguments
are needed for the cases where $\gamma_0$ is or is not part of a periodic geodesic
of period $\le T$, where $T$ is as above.

Given $\gamma\in\Pi$ we let  $T^*\gamma\subset
T^*M$ and $S^*\gamma\subset S^*M$ be the cotangent and unit cotangent
bundles over $\gamma$, respectively.  Thus, if $(x,\xi)\in T^*\gamma$ then
$\xi_\sharp$ is a tangent vector  to $\gamma$ at $x$ if $T^*M\ni\xi\to \xi_\sharp\in TM$ is the
standard musical isomorphism, which, in local coordinates, sends $\xi=(\xi_1,\xi_2)
\in T^*_xM$ to $\xi_\sharp=(\xi_\sharp^1,\xi_\sharp^2)$ with
$\xi_\sharp^j=\sum_k g^{jk}(x)\xi_k$.  Note that if $\gamma\in \Pi_{cl}$ then $L(x,\xi)\equiv
t(\gamma)<\infty$ for $(x,\xi)\in S^*\gamma$.  With this in mind, we shall let $\Pi_{cl}(T)$ denote
those $\gamma\in \Pi_{cl}$ for which $L(x,\xi)\le T$ if $(x,\xi)\in S^*\gamma$.

Let us first see that a stronger version of \eqref{2.8} must be valid whenever
$\gamma\in \Pi_{cl}(T)$.  We first note that if $A\in \Psi^0_{cl}(M)$ then
\begin{equation}\label{2.10}
A\chi(T(P-\lambda))f(x)=T^{-1}\int \Hat \chi(t/T) \, e^{-i\lambda t}
\Bigl(\, Ae^{itP}f\Bigr)(x)\, dt,
\end{equation}
and recall that because of the support properties of $\Hat \chi(t)$, the integral vanishes
when $|t|\ge T/2$.
The operator
$$f\to \Bigl(Ae^{itP}f\Bigr)(x)$$
is a Fourier integral operator with
wave front set
\begin{equation}\label{2.100}
\bigl\{\, (x,t,\xi,\tau, y,-\eta): \, \Phi_t(x,\xi)=(y,\eta), \, \pm \tau =p(x,\xi), \,
(x,\xi)\in \text{supp } A(x,\xi)\, \bigr\},
\end{equation}
where $p(x,\xi)$ is the principal symbol of $P=\sqrt{-\Delta_g}$ and $A(x,\xi)$ is the
symbol of $A$.  If $R_\gamma$ denotes the restriction to $\gamma$ then we are really concerned
with the operator
\begin{equation}\label{2.11}
f\to R_\gamma Ae^{itP}f.
\end{equation}
Regarded as an operator from $C^\infty(M)\to C^\infty
(\gamma_0\times [-T/2,T/2])$, if $\text{supp }A(x,\xi)\cap S^*\gamma = \emptyset$,
 this is a Fourier integral operator of order zero which is locally a canonical graph
\footnote{Since, for fixed $t$,  $e^{itP}: C^\infty(M)\to C^\infty(M)$ is a nondegenerate Fourier integral
operator, one needs only to check this assertion for $t=0$, in which case it is an easy
calculation using any parametrix for the half-wave operator.}.
If $\gamma=\gamma_0\in \Pi_{cl}(T)$
and we take $A=B$, where $B$ is as above, then this is automatically the case
since $B=Id-b$ and $b$ has a symbol which equals one in a neighborhood of $S^*\gamma_0$
if $\gamma_0\in \Pi_{cl}(T)$.  Therefore, by H\"ormander's \cite{hormander2}
$L^2$-estimates for nondegenerate Fourier integral operators we have
$$\int_{-T/4}^{T/4}\int_{\gamma_0}|Be^{itP}f|^2\, ds dt \le C\|f\|_{L^2(M)}^2.$$
The constant $C$ here of course depends on $T$ and $\gamma_0$ (with its main
dependence being on $\text{dist}(S^*\gamma_0, \, \text{supp }B(x,\xi))$.
Since the Fourier integral \eqref{2.11} with $A=B$
will also be nondegenerate if $\gamma$ is close to $\gamma_0$, we
conclude that whenever $\gamma_0\in \Pi_{cl}(T)$, there must be a
neighborhood ${\mathcal N}(\gamma_0)$ in $\Pi$ and a constant
$C_{\gamma_0,B,T}$ so that
$$\int_{-T/4}^{T/4}\int_\gamma |Be^{itP}f|^2\, ds dt
\le C_{\gamma_0,B,T}\|f\|_{L^2(M)}^2, \quad \gamma\in {\mathcal N}(\gamma_0).$$
If we use the Schwarz inequality and \eqref{2.10} we conclude from this that
\begin{equation}\label{2.111}
\int_\gamma |B\chi(T(P-\lambda))f|^2 \, ds \le TC_{\gamma_0,B,T}\|f\|_{L^2(M)}^2,
\quad \gamma\in {\mathcal N}(\gamma_0),
\end{equation}
which is stronger than \eqref{2.9} for these $\gamma$.

Let us now see that we also have favorable bounds on $\Pi\backslash \Pi_{cl}(T)$.
If we fix a $\gamma_0$ in this set and choose a $C\in \Psi^0_{cl}(M)$ whose symbol
vanishes on a conic neighborhood of $T^*\gamma_0$ then by the above arguments
there must be a conic neighborhood of $\gamma_0$ on which we have the analog of
\eqref{2.111} when $B$ is replaced by $C\circ B$.
This fact is independent of whether or not $\gamma_0$ is periodic.  It is just our earlier
observation that \eqref{2.11} is a nondegenerate Fourier integral operator when the symbol
of $A$ vanishes in a conic neighborhood of $T^*\gamma_0$.

 Thus, in order to show that we have uniform bounds as in \eqref{2.9} on a neighborhood of such a $\gamma_0\in \Pi\backslash \Pi_{cl}(T)$, it is
enough to show that if $A\in \Psi^0_{cl}(M)$ has a symbol supported in a small
neighborhood of $T^*\gamma_0$ then we have
\begin{equation}\label{2.12}
\int_\gamma |A\chi(T(P-\lambda))f|^2 \, ds \le CT^{-1}\lambda^{1/2}\|f\|_{L^2(M)}^2
+C'_{T,A,\gamma_0}\|f\|_{L^2(M)}^2
\end{equation}
for every $\gamma\in \Pi$.

Note that for every $x\in \gamma_0$, $T_x\gamma_0$ is one-dimensional
and if $\xi\in T_x^*\gamma$ then $-\xi\in T_x^*\gamma_0$, since $\pm \xi_\sharp
\in T_x\gamma$ are the corresponding tangent vectors to $\gamma_0$ at
$x$ pointing in opposite directions.  Thus, $T^*\gamma_0$ naturally splits
into two components, which we shall denote by $T^*_\pm \gamma_0$, and in order
to prove \eqref{2.12}, it suffices to show that the estimate holds if the symbol of $A$ is supported in
a small neighborhood of one of them, say, $T^*_+\gamma_0$, since the same argument
will apply to $T^*_-\gamma_0$.

We shall assume in what follows that the injectivity radius of $(M,g)$ is $10$ or more.
If not than we can subdivide $\gamma$ into a finite number of segments of length
smaller than one tenth of the injectivity radius and use the argument that follows to
prove the analog of \eqref{2.12} for each of these, which in turn yields \eqref{2.12}
for all of $\gamma$.

Let $Sf=A\chi(T(P-\lambda))f|_{\gamma}$ then we wish to show that 
$$\bigl(\|S\|_{L^2(M)\to L^2(\gamma)}\bigr)^2\le CT^{-1}\lambda^{1/2}+C_{T,A,\gamma_0}.$$
This is equivalent to saying that the dual operator $S^*: L^2(\gamma)\to L^2(M)$ with
the same norm, and since
$$\|S^*g\|^2_{L^2(M)}=\int_\gamma SS^*g\, \overline{g} \, ds
\le \|SS^*g\|_{L^2(\gamma)}\|g\|_{L^2(\gamma)},$$
we would be done if we could show that 
\begin{equation}\label{2.13}
\|SS^*g\|_{L^2(\gamma)}\le
\Bigl(\, CT^{-1}\lambda^{1/2}+C_{T,A,\gamma_0}\, \Bigr)\, \|g\|_{L^2(\gamma)}.
\end{equation}
But the kernel of $SS^*$ is $K(\gamma(s),\gamma(s'))$
where $\gamma(s)$ parameterizes $\gamma$ by arclength and $K(x,y)$,
$x,y\in M$ is the kernel of the operator $A\circ \rho(T(P-\lambda))\circ A^*$
with $\rho = (\chi(\tau))^2$ being the square of $\chi$.  Its Fourier transform
$\Hat \rho$ is the convolution of $\Hat \chi$ with itself and thus
$\Hat \rho(t)=0$, $|t|\ge1/2$.  Consequently, we can write
\begin{equation}\label{2.14}
A\circ \rho(T(P-\lambda))\circ A^* = T^{-1}\int_{-T/2}^{T/2}
\Hat \rho(t/T) \, e^{-it\l} \, \Bigl(A\circ e^{itP}\circ A^*\Bigr)\, dt.
\end{equation}

The wave front set of the kernel of 
$$A\circ e^{itP}\circ A^*$$
regarded as an operator from $C^\infty(M)$ to $C^\infty(M\times {\mathbb R})$ is contained in
\begin{equation}\label{WF}
\{(x,t,\xi,\tau; y,-\eta): \, \Phi_t(x,\xi)=(y,\eta), \tau =\pm p(y,\eta), \, 
(x,\xi), (y,\eta)\in \text{supp }A\}.
\end{equation}
Our assumption that $\gamma_0\notin \Pi_{cl}(T)$ 
implies that if $(x,\xi)\in S^*\gamma_0$
then $\{\Phi_t(x,\xi): \, 1\le |t|\le T\}$ must be a closed subset of $S^*M$ which is disjoint
from $\{(x,\xi)\}$.  If we assume also that $(x,\xi)$ and $(y,\eta)$ belong to the
same component $S^*_+\gamma_0$ of $S^*\gamma_0$ then we have that $\Phi_{t_0}(x,\xi)
=(y,\eta)$ for some $|t_0|\le1$, and therefore $\Phi_t(x,\xi)\ne (y,\eta)$,
 for $|t|\in [2,T-1]$ since if $\Phi_{t}(x,\xi)=(y,\eta)$ then we
must also have $\Phi_{t-t_0}(x,\xi)=(x,\xi)$ for this $t_0$.  
Consequently,
$$\{(x,\xi, \Phi_t(x,\xi)): \, (x,\xi)\in S^*_+\gamma_0, \, 2\le |t|\le T-1\}
\cap S^*_+\gamma_0 \times S^*_+\gamma_0 = \emptyset,$$
and since both are compact subsets of $S^*M\times S^*M$, we deduce from \eqref{WF}
that if the symbol
of $A$ is supported in a small conic neighborhood of $S^*_+\gamma_0$,
then $K(t,x,y)$ will be $C^\infty$ when $|t|\in [2,T-1]$.

Therefore, for such $A$, if 
if $\alpha\in
C^\infty_0({\mathbb R})$ equals one if $|t|\le 3$ and zero for
$|t|\ge4$, the difference between the kernel $K(x,y)$ in \eqref{2.14} and
$$K_0(x,y)=T^{-1}\int \alpha(t) \Hat \rho(t/T) e^{-it\l}\, \Bigl(
A\circ e^{itP}\circ A^*\Bigr)(x,y)\, dt,
$$
must be bounded, by a constant 
which is independent of $x$ and $y$ (but depends on $T$, $\gamma_0$ and $A$).
 Since we are assuming that the injectivity radius of
$(M,g)$ is $10$ or more one can use the Hadamard parametrix
construction for the wave equation and standard stationary phase
arguments (cf. Chapter~5 in \cite{soggebook} or the proof of
Lemma~4.1 in \cite{burq}) to see that for $x,y\in M$ we have
$$|K_0(x,y)|\le CT^{-1}\lambda^{1/2}\, (\, \text{dist}_g(x,y)\, )^{-1/2}+C_A.$$
Since this kernel restricted
to $\gamma\times \gamma$ gives rise to an integral operator satisfying the estimates
in \eqref{2.13}, we conclude that we also have uniform bounds of the
form \eqref{2.9}, when $A$ is as above.

This completes the proof that the analog of \eqref{2.9} holds for all
$\gamma$ in some neighborhood of $\gamma_0$ when $\gamma_0\in
\Pi \backslash \Pi_{cl}(T)$.

Combining what we have done for $\Pi_{cl}(T)$ and $\Pi\backslash
\Pi_{cl}(T)$, since $\Pi$ is compact, we conclude that \eqref{2.9}
must be valid with uniform constants for every $\gamma\in \Pi$.
This completes the proof of \eqref{2.8} and hence \eqref{2.7}.
Since the latter holds for all $\varepsilon>0$, we conclude from
\eqref{2.3} that
\begin{equation}\label{2.15}
\limsup_{\lambda\to\infty}\frac1{N(\lambda)}\sum_{\l_j\le \l}\bigl(\,
\l_j^{-1/8}\|e_{\l_j}\|_{L^4(M)}\, \bigr)^{12}=0.
\end{equation}

We can now finish the proof of Theorem~\ref{theorem1.1} using a counting
argument from \cite{cdv} and \cite{zelditch}.  If $S\subset\{ \lambda_j \}$,
we define its density to be
$$D(S)=\liminf_{\l\to \infty}\frac{\#\{\l_j\in S: \, \l_j\le \l\}}{N(\l)}.$$
Then if we use \eqref{2.15} we conclude that for every $n=2,3,\dots$ we can
find a subset $S_n$ of the eigenvalues $\{\l_j\}$ so that
$$D(S_n)\ge 1-\frac1n \, \, \, \text{and } \, \,
\l_j^{-1/8}\|e_{\l_j}\|_{L^4(M)}\le \frac1n, \, \, \l_j\in S_n.$$
Using this we conclude that there must be a set $S_\infty
=\{\l_{j_k}\}\subset \{\l_j\}$ of density $1$ so that
$$\limsup_{k\to\infty}\l_{j_k}^{-1/8}\|e_{\l_{j_k}}\|_{L^4(M)}=0.$$
Indeed, by the above, we can choose increasing $N_\nu\in {\mathbb N}$,
$\nu=2,3,\dots$ so that
$$ \#\{\l_j\in S_\nu: \, \l_j\le \l\}/N(\l)\ge 1-2/\nu, \quad \forall \, \l\ge N_{n-1}.$$
Consequently,
$$S_\infty =\bigcup_{\nu=2}^\infty S_\nu \cap \{\l_j \le N_\nu\}$$
will have the desired properties. \qed


\newsection{Average $L^4$ norms of spherical harmonics}

In this section we shall prove \eqref{1.7}:

\begin{theorem}\label{theorem2}  Let $\{Y^k_m(x)\}_{m=-k}^k$,
$x=(\sin\phi \cos\theta, \sin \phi \sin\theta, \cos\phi)$, be the
orthonormal
basis of spherical harmonics of degree $k$ defined in \eqref{1.6}.  Then  there is a uniform constant $C$
so that
\begin{equation}\label{3.1}
\frac1{2k+1}\sum_{m=-k}^k\int_{S^2}|Y^k_m|^4dV  \le C\log k, \quad
k\ge2.
\end{equation}
Moreover, if $\1=(0,0,1)$ and if $r=\min_\pm \text{dist}(x,\pm\1)$
\begin{equation}\label{3.2}
\Bigl(\sum_{m=-k}^k |Y^k_m(x)|^4\Bigr)^{1/4} \le \begin{cases}
Ck^{1/4} r^{-1/4}\bigl(\log (kr)\bigr)^{1/4}, \quad r\ge 2/k
\\ \\
Ck^{1/2}, \quad r\le 2/k.
\end{cases}
\end{equation}
\end{theorem}

Clearly \eqref{3.2} implies \eqref{3.1}, and so we just need to prove the second inequality
in the theorem.  To prove this, we first realize that by Parseval's theorem we have
$$2\pi \sum_{m=-k}^k |Y^k_m(x)|^4 =\int_0^{2\pi} \Bigl| \, \sum_{m=-k}^k |Y^k_m(x)|^2
e^{im\theta}\Bigr|^2\, d\theta.$$
The kernel $\Pi_k(x,y)$ for projection on to spherical harmonics of degree $k$ is given
my the formula
$$\Pi_k(x,y)=\sum_{m=-k}^k Y^k_m(x)\overline{Y^k_m(y)},$$
which means that
\begin{equation}\label{3.3}
2\pi \sum_{m=-k}^k |Y^k_m(x)|^4 =\int_0^{2\pi}
|\Pi_k(x,e^{i\theta}x)|^2 \, d\theta,
\end{equation}
if we abuse notation a bit and let $e^{i\theta}x$ denote rotation of our
vector $x=(x_1,x_2,x_3)$ by angle $\theta$ about the $x_3$-axis,
i.e., $e^{i\theta}x= (\cos\theta x_1,\sin\theta x_2,x_3)$.
Using the well known bounds (see \cite{szego}, \cite{sph}) for $\Pi_k$,
$$|\Pi_k(x,y)|\le Ck^{1/2}(k^{-1}+\text{dist}(x,y))^{-1/2}, \quad x,y\in S^2$$
we conclude that
\begin{equation}\label{3.4}\sum_{m=-k}^k|Y^k_m(x)|^4 \le Ck
\int_0^{2\pi} \bigl(k^{-1}+\text{dist}(x,e^{i\theta}x))^{-1} \, d\theta.
\end{equation}
Since $\text{dist}(x,e^{i\theta})\le Cr |\sin\theta|$, we conclude that the
right side of \eqref{3.4} is $\le Ckr^{-1}\log(kr)$ if $r\ge2/k$, and $\le Ck$
if $r\le 2/k$,
for some uniform constant $C$ when $k\ge2$, which is just \eqref{3.2}.  \qed

We believe this estimate to be sharp, but defer the analysis to
the future.

\subsection{\label{ONB} $L^4$ norms of orthonormal bases of spherical harmonics}

The space of  ${\mathcal ONB}_k$ of  Hermitian  orthonormal bases
of ${\mathcal H}_k$ may be identified with the unitary group $U(2k
+ 1)$. Any orthonormal basis $\Phi_k = \{\phi^k_1, \dots,
\phi^k_{2k + 1}\}$ can be obtained by applying an element $U \in
U(2k + 1)$ to the standard orthonormal basis $\{Y^k_m\}$. We then
consider the functional on ${\mathcal ONB}_k$ defined by,
$$\Lambda_k^4(\Phi_k) = \sum_{m =1}^{2k +1} ||\phi^k_{ m}||_{L^4}^4. $$
We have just proved that $\Lambda_k^4$ of the standard orthonormal
basis is bounded by $k \log k$.

One may consider similar functionals on orthonormal bases of all
eigenspaces with $\lambda_k \leq \lambda$, i.e. the direct sum
$\bigoplus_{k \leq \lambda} {\mathcal H}_{k}$. The functional then
has a natural generalization to any $(M, g)$ and is essentially
the one studied in previous sections.

\subsection{\label{RANDOM} Random orthonormal bases of spherical harmonics}

We now consider the $\Lambda_k^4$ functional on a random basis of
spherical harmonics. The question we pose is, what is the average
value of the functional on random orthonormal bases? In \cite{SZ}
we considered problems of this kind for $L^p$ norms of individual
eigenfunctions, but there is a new dimension to the problem for
random orthonormal bases. For background on random orthonormal
bases we refer to \cite{SZ}.

 We introduce the probability space
$({\mathcal ONB},d\nu)$, where ${\mathcal ONB}$ is the infinite
product of the sets, and $\nu =\prod_{k=1}^{\infty} \nu_k$, where
$\nu_k$ is Haar probability measure on ${\mathcal ONB}_k$. A point
of ${\mathcal ONB}$ is thus a sequence ${\bf \Phi} = \{(\phi^k_1,
\dots, \phi^k_{2k + 1}) \}_{k\geq 1}$ of orthonormal bases of
${\mathcal H}_k$.

The functionals we are interested in are
\begin{equation}\label{A} \Lambda_{k}^4({\bf \Phi}) = \sum_{j = 1}^{2k+ 1} \int_{S^2} \|\phi^k_j(x))\|^4 dV\;.\end{equation}
If we fix the standard ONB $Y_k =\{Y^k_{m}\}$ and express every
other as $U_k Y^k$, then our functional is
\begin{equation}\label{AU} \Lambda_{k}^4(U) = \sum_{j = 1}^{2k + 1} \int_{S^2} \| (UY^k)_j(z)\|^4 dV.\end{equation}

Let $d\mu_k$ be normalized Haar measure on $U(k)$ and let $\E_k$
denote expectation with respect to this measure. We conjecture
that \begin{equation} \label{CONJ} \E_k \Lambda_k^4 =  2 (2 k +
1),
\end{equation}
i.e. the elements on average have $L^4$ norm equal to $2$.

We briefly sketch the proof.  We start from the fact that
\begin{equation}\label{AUa} \begin{array}{lll} \E \Lambda_{k}^4(U) &= & \sum_{j = 1}^{2k + 1} \sum_{m_1, m_2, m_3, m_4 = -
k}^k \left(\int_{S^2} Y^k_{m_1}(x) \overline{Y}^k_{m_2}(x)
Y^k_{m_3}(x) \overline{Y}^k_{m_4}(x) dV \right)\\ && \\ &&
\left(\int_{U(2k + 1)} U^j_{m_1} \overline{U}^j_{m_2} U^j_{m_3}
\overline{U}^j_{m_4} d\mu_k(U) \right). \end{array}
\end{equation}
In fact, the sum over $j$ is constant, so the right side equals
$(2k + 1)$ times
\begin{multline*} \sum_{m_1, m_2, m_3, m_4 = -
k}^k \left(\int_{S^2} Y^k_{m_1}(x) \overline{Y}^k_{m_2}(x)
Y^k_{m_3}(x) \overline{Y}^k_{m_4}(x) dV \right)
\\
 \times \left(\int_{U(2k +
1)} U^1_{m_1} \overline{U}^1_{m_2} U^1_{m_3} \overline{U}^1_{m_4}
d\mu_k(U) \right).  \end{multline*}

The integrals $\left(\int_{U(2k + 1)} U^j_{m_1}
\overline{U}^j_{m_2} U^j_{m_3} \overline{U}^j_{m_4} d\mu_k(U)
\right)$ were first  studied by Weingarten  \cite{W}. The main
result is that the random variables $\{\sqrt{2k + 1} U_{i j}\}$
behave asymptotically like independent complex Gaussian random
variables of mean zero and variance one. Exact formulae are given
in \cite{CS}, and the latter can be used to determine the
asymptotics of our sums over $2k + 1$ indices with different
coefficients as $k \to \infty$. The dominant terms come from the
cases where all $m_j$ are equal (then one has the fourth moment of
the Gaussian) or when the indices $m_j$ are paired into couples
(one barred and one unbarred). Then we have,
\begin{equation}\label{AUb} \begin{array}{l}(2k + 1)^2 \E \Lambda_{Nk}^4(U) \\ \\  \sim  (2 k+ 1)\left( 2 \sum_{m
= -k}^k \int_{S^2} |Y^k_{m}(x)|^4 dV  + 2 \sum_{m_1 \not= m_2 = -
k}^k \left(\int_{S^2} |Y^k_{m_1}(x)|^2 |Y_{m_2}^k|^2 dV \right) \right) \\  \\
= 2 (2k + 1) \int_{S^2} |\Pi_k(x,x)|^2 dV = 2 (2k + 1)^3.
\end{array}
\end{equation}
Here, we use that the 4th moment of the complex normal Gaussian
equals $2$ and that there are two ways to pair the indices in the
off diagonal terms. Dividing by $(2k + 1)^2$ then implies the
result.

\subsection{Other orthonormal bases}

Theorem \ref{theorem2} shows that the $\Lambda_k^4$ functional on
the standard basis $\{Y^k_m\}$ is only $\log k$ higher than for a
random orthonormal basis.  Hence it is doubtful that it does not
maximize $\Lambda_k^4$. We do not know which ONB maximizes the
functional, but in this section we suggest a possible construction
of one which has a higher $\Lambda_k^4$ value than the standard
basis.

As mentioned above,  the highest wight spherical harmonic $Y^k_k $
has $L_4^4$ equal to $k^{1/2}$ on $S^2$, thus maximizing  the norm
functional. This suggests constructing orthonormal bases
$\phi^k_{\gamma_j}$ consisting in part  of highest weight
spherical harmonics with respect to a well-separated set of closed
geodesics $\gamma_j$. That is, for each closed geodesic $\gamma$,
one introduces the subgroup $G_{\gamma}$ of rotations fixing
$\gamma$ (as a set) and then constructs $Y^k_m$'s with respect to
this circle action.

Of course, the $\phi^k_{\gamma_j}$ are not orthogonal, and their
inner products $\langle \phi^k_{\gamma_j}, \phi^k_{\gamma_i}
\rangle$ depend on the angle $\vartheta_{j, k}$ between the
geodesics $\gamma_j, \gamma_i$. To construct an orthonormal basis
it would be necessary to apply Gram-Schmidt to such
$\phi_{\gamma_j}^k$, and in the process one may destroy the high
$L^4$ norms of the resulting eigenfunctions. The question is, how
many $\phi^k_{\gamma_j}$ can be used in such a construction while
preserving the high $L^4$ norms of these Gaussian beams?

The geodesics $\gamma_j$ are points in the space $G(S^2, g_0)$ of
geodesics of $S^2$. A well-separated set of $2k + 1$ geodesics
(i.e. a basis) would only have separation of order
$k^{-\frac{1}{2}}$. To beat the bound $k \log k $ for the standard
basis one it would suffice to construct a partial  orthonormal
basis containing $k^{1- \delta} $ roughly Gaussian beams  with
roughly $||\phi^k_{\gamma_j}||_4^4 \simeq \sqrt{k}$ and with
$\delta < \frac{1}{2}$. One would then complete it with an
arbitrary orthonormal basis of the ortho-complement of the span.
It would be interesting to see how far separated the $\gamma_j$
would need to be so that Gram-Schmidt would not destroy the bounds
$||\phi^k_{\gamma_j}||_4^4 \simeq \sqrt{k}$ too much.


\end{document}